 
 
 
 

\documentstyle{amsppt}
\magnification 1200

\newcount\numofeq 
\newcount\numofsec 
\newcount\numofthrm 
\numofsec=0
\numofeq=0
\numofthrm=0

\define\equano{\global\advance\numofeq by 1\ifnum\numofsec >0 
\eqno(\the\numofsec.\the\numofeq) \else\eqno(\the\numofeq) \fi}
\define\numnow{\ifnum\numofsec >0 ({\the\numofsec}.{\the\numofeq})\else
({\the\numofeq})\fi}
\define\numbefore{\advance\numofeq by -1 \ifnum\numofsec >0
({\the\numofsec}.{\the\numofeq})\else({\the\numofeq})\fi}
\define\numbeforethat{\advance\numofeq by -2 \ifnum\numofsec >0 
({\the\numofsec}.{\the\numofeq})\else({\the\numofeq})\fi}
\define\numback#1{\advance\numofeq by 
-#1\ifnum\numofsec>0({\the\numofsec}.{\the\numofeq})\else({\the\numofeq})\fi}

\define\R{\Bbb R}
\define\norm{|\!|}
\define\norma#1{\norm #1\norm}
\define\({\left (}
\define\){\right )}
\define\innerprod#1#2{\langle#1,#2\rangle}  
\global\numofthrm=1
\global\numofeq=0
\define\pr{{r\over\sqrt{1+r^2}}}
\define\8{\infty}

\topmatter

\title On zonoids whose polars are zonoids \endtitle

\author Yossi Lonke \endauthor
\leftheadtext{Y. Lonke}
\address Institute of mathematics, Hebrew University, Jerusalem 91904, Israel 
\endaddress
\email lini\@math.huji.ac.il \endemail
\thanks Supported in part by the United States-Israel Binational
 Science Foundation \endthanks
\abstract
 Zonoids whose polars are zonoids, cannot have proper faces
 other than vertices or facets. However, there exist non--smooth
zonoids whose polars are zonoids. Examples in $\R^3$ and $\R^4$ are given.
\endabstract
\endtopmatter
\document
\openup 2pt
\tolerance=600

\head  Introduction \endhead
A {\it zonotope\/} in $\R^n$ is a vector sum of segments. A {\it zonoid\/} in 
$\R^n$ is a limit of zonotopes in $\R^n$
with respect to the Hausdorff metric. The sum of the centers of the segments, 
defines
a center of symmetry for a zonotope, and so by definition every zonoid is
centrally symmetric, compact and convex. Consequently, every zonoid is the unit
ball of some norm. The special structure of zonoids allows a more precise 
statement about what kind of norms enter the discussion of zonoids. Writing down
the support function of a zonotope we see that zonoids are
precisely unit balls of quotients of $L_{\8}$ spaces.
Since every two dimensional convex, centrally symmetric and compact
set is a zonoid, every two dimensional Banach space is isometric to a
subspace of $L_1$. Therefore there is no zonoid theory in $\R^2$. For detailed 
discussions concerning zonoid theory see \cite{1,5,12}.\par
A well known theorem in Functional Analysis, due to Grothendieck, asserts that
among infinite dimensional Banach spaces, the ones which are isomorphic both
to a subspace of $L_1$ and to a quotient--space of $L_{\8}$, are isomorphic
to a Hilbert space. (For a proof see \cite{8}). A natural question is whether 
this
"isomorphic" theorem
 has an "isometric" analogue. E. Bolker conjectured in \cite{1} (conjecture 6.8) 
that there {\it is\/} a finite
dimensional isometric version of Grothendieck's theorem. Bolker formulated his
conjecture in the language of zonoids, which amounts to the conjecture
that a zonoid whose dimension is at least~$3$ and whose polar is also a zonoid 
must be an ellipsoid.\par
Six years later, R. Schneider in \cite{10}, constructed examples of zonoids 
whose polars are also zonoids,
which are not ellipsoids, and consequently Bolker's conjecture was disproved.
The infinite dimensional isometric problem is still open.\par
Schneider used spherical harmonics to prove that it is possible to apply
smooth perturbations to the Euclidean ball in $\R^n$, such that the resulting 
bodies
are zonoids whose polars are zonoids. Since the set of zonoids in $\R^n$, where 
$n\geq  3$, is closed and nowhere
dense with respect to the Hausdorff metric, the perturbations had to be 
performed with
respect to an essentially different metric. The metric employed by Schneider 
involved
high--order derivatives of the support function. It seemed plausible after 
Schneider's
work that zonoids whose polars are zonoids must be smooth, and hence strictly 
convex.\par
In \S2 of this work an example of a non--smooth zonoid
whose polar is also a zonoid is presented. It consists of forming the
Minkowski sum of the~$3$--dimensional Euclidean unit ball and a
concentric circle of radius~$1$. Similar examples exist in $\R^4$. 
\par
Although smoothness of zonoids whose polars are zonoids
cannot be guaranteed, it nevertheless cannot be lost in an arbitrary fashion.
\S1 contains some information in this direction.
The result is

\proclaim{Theorem} Suppose $K=B+C$ is a convex body in $\R^n$, where $n\geq 3$
and $B,C$ are convex, compact and
centrally symmetric subsets of $\R^n$. If  $1\leq{\text {dim\,}} C\leq n-2$
then the polar of $K$ is not a zonoid.\endproclaim

It is well known that every face of a zonoid $Z$ is a translate of a zonoid 
which
is a summand of $Z$. (\cite{1}, Th. 3.2, and also \cite{11}, p.189). Therefore, 
the theorem implies:

\proclaim{Corollary} If $n\geq 3$ and $Z$ is an $n$--dimensional zonoid whose
polar is also a zonoid, then the boundary of~$Z$ does not contain proper faces 
whose
dimension is different from $n-1$ or zero.\endproclaim

As far as the dimension is concerned, the examples of \S2 show that these
statements cannot be improved.\par
Another immediate corollary of the theorem is that the polar of a zonotope whose
dimension is at least~$3$ is not a zonotope. This had been proved long ago
by M.A. Perles, and by E. Bolker. (cf. \cite{1}). Both proofs are based on the 
special
polytopal structure of a zonotope. In particular, Perles shows that every 
zonotope
whose dimension is at least~$3$ has more vertices than facets.\par

\head \S1. Proof of the theorem \endhead
\advance\numofsec by 1
A fundamental property of zonoids which is exploited below
appears as  Theorem 3.2 in \cite{1}. It is stated here as a lemma. 

\proclaim{Lemma 1.1} Every proper face of a zonoid $K$ is a translate of a 
zonoid of lower
dimension which is a summand of~$K$.\endproclaim

A summand can be either direct or not. $B$ is said to be a direct summand of $K$
if ${K=B+C}$ and ${{\text {dim\,}} B+{\text {dim\,}} C={\text {dim\,}} K}$. In 
such cases the Minkowski sum is written in the form 
${K=B\oplus C}$. Due to this distinction between types of summands, the proof
of the theorem will be divided into two parts. The first part consists of a 
proposition
which settles the case of direct summands. The restriction ${\text {dim\,}} 
C\leq n-2$ which
appears in the formulation of the theorem does not play any role in the setting
of direct summands. 

\proclaim{Proposition 1.2} A polar of a zonoid whose dimension is at least~$3$ 
does not
have non--trivial direct summands.\endproclaim

\demo{Proof}
Assume that $K$ is a polar of
a zonoid, ${\text {dim}}\,K\geq 3$ and $K$ does have non--trivial direct 
summands.
Applying a suitable linear transformation, it can be assumed that $K=B\oplus C$
 where ${\text {span}}\,B$ and ${\text {span}}\,C$ are mutually orthogonal.
Let~$P$ denote the orthogonal projection onto ${\text {span}}\,B$
and let~$Q=I-P$. Then the norm of~$K$ can be written as 
$$\norma{x}_K=\max\{\norma{Px}_B,\norma{Qx}_C\},\enspace\forall 
x\in\R^n.\equano$$
Choose any point $x$ on the boundary of~$B$, and any point~$y$ on the boundary
of~$C$. Then 
$$\norma{(x+y)+(x-y)}_K+\norma{(x+y)-(x-y)}_K=4=2\(\norma{x+y}_K+\norma{x-y}_K\)
.\equano$$
By assumption, $K^{\circ}$ is a zonoid. Hence there is a positive measure $\nu$
 on the sphere such that 
$$\norma{z}_K=\int_{S^{n-1}}|\innerprod {z}u|\,d\nu(u),\enspace\forall 
z\in\R^n.\equano$$
For every $x\in\partial B$ and $y\in\partial C$ consider two functions defined
on the sphere by
$$f_{y,x}(u)=\innerprod {x+y}u\enspace\hbox{and}\enspace g_{y,x}(u)=\innerprod 
{x-y}u,$$
Then by $\numnow$ and $\numbefore$, for
the norm in $L^1(d\nu)$, and for $f=f_{x,y}$ and $g=g_{x,y}$, one has
$$\norma{f+g}+\norma{f-g}=2\(\norma{f}+\norma{g}\).\equano$$
Such an equality can occur only if
$$f(u)g(u)=0\quad\hbox{for $\nu$--almost every $u\in S^{n-1}$}.\equano$$
This implies: 
$${\text {supp}}\,\nu\subset\{u\in S^{n-1}:|\innerprod ux|=|\innerprod 
uy|\}\quad 
\forall\,y\in\partial C,\forall\,x\in\partial B.\equano$$
Observe that if $x_1,\dots,x_n$ is a linear basis for~$\R^n$ then the set
$$\{u\in S^{n-1}:|\innerprod u{x_1}|=|\innerprod u{x_2}|=\cdots=|\innerprod 
u{x_n}|\}\equano$$
is finite (it contains at most $2^n$ points). Since the dimension of $K$ is at 
least three, the dimension of one of the 
summands must be at least two, and so $\numbefore$ cannot hold for every choice
of points $x\in\partial B$, and $y\in\partial C$, due to the previous
observation. This contradiction proves the proposition.\qed\enddemo

The second part of the proof consists of proving the theorem for non--direct
summands. Our original proof was for a $1$--dimensional summand, i.e, a segment,
but as was kindly pointed out by R. Schneider, the same argument yields
the result as stated here.\par
Let $K$ be a centrally symmetric, compact and convex subset of $\R^n$.
Given any subset $U\subset\R^n$,
consider the following subset of the boundary of~$K\,$:
$$A(K,U)=\{x\in\partial K:\hbox{$K$ has an outer normal at~$x$ which is
orthogonal to {\text {span}}\thinspace $U$}\}.$$
Concerning such sets the following Lemma holds.

\proclaim{Lemma 1.3} Suppose $K$ is a centrally symmetric compact and convex
subset of $\R^n$, such that ${\text {\rm dim\,}} K\geq 2$.
Let $U\subset \R^n$ denote a subset whose span is of dimension at most $n-2$.
Then for every $(n-1)$--dimensional subspace~$H$ 
 the intersection ${A(K,U)\cap H}$ is not empty.\endproclaim

\demo{Proof} 
One may assume that $U$ is a subspace of dimension $n-2$. Let $H$ be a subspace 
of dimension $n-1$.\par
First assume that $K$ is smooth and strictly convex. Let $G_K$ denote the Gauss 
map, which takes every point $x$ on $\partial K$ to the outward unit normal 
of~$K$ at~$x$. Then $G_K$ is a homeomorphism between $\partial K$ and 
$S^{n-1}$.\par
The intersection $U^{\perp}\cap S^{n-1}$ is a great circle. Since the set 
$A(K,U)$ is the inverse image of this intersection under the Gauss map, it is a 
connected, centrally symmetric subset of the boundary of $K$ and hence meets 
$H$.
This proves the assertion under the special assumption on~$K$.\par
The proof for general $K$ is now concluded by an approximation argument. Let $K$ 
be as in the formulation of the lemma. Choose a sequence $(K_i)_{i\in\Bbb N}$ of 
centrally symmetric, smooth and strictly convex bodies converging to~$K$. For 
each $i\in\Bbb 
N$ there exists a pair $(x_i,u_i)\in H\times U^{\perp}$ such that 
$x_i\in\partial K_i$ and $u_i$ is an outer normal vector $K_i$ at $x_i$,
thus $h(K_i,u_i)=\innerprod {x_i}{u_i}$, where $h$ denotes the support function. 
The sequence $(x_i,u_i)$ has a convergent subsequence, and one may assume that 
this sequence itself converges to a pair $(x,u)$. Since the support function $h$ 
is simultaneou
sly continuous in both variables, one gets $h(K,u)=\innerprod xu$. If 
$h(K,u)\neq 0$, then $x\in\partial K$, and from $(x,u)\in H\times U^{\perp}$ one 
gets $x\in A(K,U)\cap H$. On the other hand, if $h(K,u)=0$, then
since ${\text {dim\,}} K\geq 2$, there exists a point $y\in\partial K\cap H$, 
and this point is in the set $A(K,U)\cap H$. \qed\enddemo

\medskip\noindent
{\bf Proof of the theorem}\par\noindent 
 The main idea is similar to the one which appeared above, in the
proof of the proposition. Suppose $K^{\circ}$ (the polar of~$K$) is a zonoid, 
$K$
is $n$--dimensional and has a non-direct summand $C$ where ${1\leq {\text 
{dim\,}} C\leq n-2}$. 
Then $K=B+C$, for some centrally symmetric compact and convex subset~$B$
of $\R^n$. For every $x\in A(B,C)$, the
set $x+C$ lies entirely on the boundary of $K$. Therefore,
$$\norma{(x+c)+(x-c)}_K=2=\norma{x+c}_K+\norma{x-c}_K,\enspace\forall
x\in A(B,C),\,\forall c\in C.\equano$$
Arguing similarly as in the proof of the proposition, an equality in the 
triangle
inequality in $L_1(S^{n-1},\nu)$ is obtained, where the vectors are the 
functions 
$$f_{x,c}(u)=\innerprod {x+c}u\enspace\hbox{and}\enspace g_{x,c}(u)=\innerprod 
{x-c}u,\equano$$
and $c\in A(B,C)$ is arbitrary. Since $K^{\circ}$ is a zonoid, there exists a 
positive Borel measure $\nu$ on
the sphere which satisfies an equation of the form $\numback 6$.
Being 
positive, the measure must assign all its mass to the set of points where
the functions from $\numnow$ have
the same sign. Therefore,
$${\text {supp}}\,\nu\subset\{u\in S^{n-1}:|\innerprod ux|\geq  |\innerprod 
uc|\}\quad\forall x\in A(B,C),\,\forall c\in C.\equano$$
By lemma 1.3, for every $u\in S^{n-1}$ there
exists a point $x\in A(B,C)$ such
that $x\perp u$, and so from $\numnow$ the measure is seen to be
concentrated on the section $C^{\perp}\cap S^{n-1}$, whose dimension is 
$(n-{\text {dim\,}} C)$. 
But this contradicts the fact that the dimension of~$K$ is~$n$. This completes
the proof of the theorem.\qed\medskip
\vfill\eject

\head \S2. The Barrel zonoid \endhead
\advance\numofsec by 1
\global\numofeq=0
The purpose of this section is to present an example of a non--smooth
zonoid whose polar is also a zonoid. 
First, the so called "barrel zonoid" is introduced
and some of its properties are discussed. Afterwards the analytic tools to be 
used are presented, followed by a calculation which yields the desired example.
The main tool here is an inversion formula for the cosine transform which
involves the Radon transform, due to Goodey and Weil \cite{4}.\medskip
Let $B_2^n$ denote the $n$--dimensional euclidean unit ball. For a positive
number $r>0$, consider the zonoid $\Cal B_{n,r}=B_2^n+rB_2^{n-1}$. It is 
invariant under rotations which keep the $n$'th coordinate fixed. Such
bodies are called {\it rotationally symmetric\/}. If $0\leq \varphi\leq \pi$ 
denotes the vertical
angle in spherical coordinates, then the restriction of the norm of $\Cal 
B_{n,r}$
to the unit sphere depends only on~$\varphi$. Therefore it can be identified 
with
a function $f_r$ defined on the interval $[0,\pi]$, and by symmetry, attention
can be restricted to the interval $[0,\pi/2]$. The rotational symmetry implies
that for all dimensions $n\geq 3$, the norm of $\Cal B_{n,r}$ is represented by 
the same function $f_r$. A simple $2$--dimensional calculation shows that 
$$f_r(\varphi)=\cases \cos\varphi,& {\text {if }} 0\leq \varphi\leq 
\tan^{-1}r\cr
{1\over r\sin\varphi+\sqrt{1-r^2\cos^2\varphi}},& {\text {if }} \tan^{-1}r\leq 
\varphi\leq \pi/2\endcases\equano$$
In case $n=3$ and $r=1$ the resulting body is barrel--shaped. Henceforth the 
name "barrel zonoid" will refer
to a body of the form $\Cal B_{n,r}$, and $\Cal B_{n,1}$ will be denoted by 
$\Cal B_n$.\par

The support function of the barrel--zonoid is the sum of the support functions
of its summands. Therefore its restriction to the sphere is the function
$1+r\sqrt{1-u_n^2}$, where $u_n$ denotes the $n$'th coordinate. 
It is not differentiable at the points $\pm e_n$
which geometrically means that there is no unique supporting hyperplane
to the polar at the points $\pm e_n$. This of course corresponds to the fact 
that
$\Cal B_{n,r}$ itself is not strictly convex. The polar $\Cal B_{n,r}^{\circ}$ 
is also rotationally symmetric,
and in case $n=3,r=1$, when it is intersected by a plane parallel
to $(0,0,1)$ and passing through the origin, the result is a symmetric, 
parabolic
curve whose equation is easily obtained
by means of the radial function of the polar, and is given by $|y|={1-x^2\over 
2}$,
for $|x|\leq 1$. Rotating this curve about the interval $|x|\leq 1$ yields the 
polar of
the $3$--dimensional barrel, (of radius $r=1$), and so an explicit figure of the 
polar may be
obtained. Its shape resembles that of an American football.\par
\medskip
 Let $\norma{\cdot}$ denote a norm in $\R^n$. In order to prove that the polar
of the unit ball determined by the given norm is a zonoid, one needs to find a
positive, symmetric measure $\mu$ on the
sphere such that
$$\norma{u}=\int_{S^{n-1}}|\innerprod uv|\,d\mu(v),\quad \forall 
u\in\R^n.\equano$$
The r.h.s of this equation is the {\it cosine transform\/} of the measure~$\mu$.
Usually
the cosine transform $T:C_e^{\8}(S^{n-1})\to C_e^{\8}(S^{n-1})$ is defined
on the space of infinitely differentiable even functions on the unit sphere
$S^{n-1}$ by
$$(Tf)(u)=\int_{S^{n-1}}|\innerprod uv|f(v)\,d\lambda_{n-1}(v),$$
where $d\lambda_{n-1}$ is the spherical Lebesgue measure on $S^{n-1}$. It is 
clear that the same
formula can be used to transform more general objects than $C^{\8}$ functions
on the sphere, such as measures. Therefore the equation $\numnow$
can be viewed as

$$T\mu = \norma{\cdot}.\equano$$
This equation is well known and has been the subject of many
investigations. In 1937, A.D. Alexandrov proved that there is at most one 
symmetric measure
which solves $\numnow$. Since then several other proofs of the same fact were 
found. See e.g.,\cite{7}.\par
There does not always exist a symmetric measure $\mu$ which solves
$\numnow$ for a given norm. In fact, it is known
that if the norm is of a polytope, then a solution exists only if this
polytope is a zonotope. (\cite{11}), corollary 3.5.6, pp. 188) However, Weil 
\cite{13} showed
that for every norm there exists a symmetric distribution~$\rho$, (i.e., a 
continuous linear
functional on the space $C_e^{\8}(S^{n-1})$) whose domain can be extended
to include the functions $|\innerprod u{\cdot}|,u\in S^{n-1}$, such that
$\rho(|\innerprod u{\cdot}|)=\norma{u}$. The fact that every distribution can
be viewed as the cosine transform of a distribution follows from the 
 self--duality of~$T:C_e^{\8}(S^{n-1})\to C_e^{\8}(S^{n-1})$ (a simple
consequence of Fubini's theorem) and a result by Schneider, asserting that
$T$ is onto $C_e^{\8}(S^{n-1})$.
Therefore $\numnow$ can always be solved with a distribution instead of a 
measure,
for any given norm, and the symbol $T^{-1}(h_K)$ acquires a precise meaning
for every given support function of a centrally symmetric convex body~$K$.
The distribution 
$T^{-1}(h_K)$ is called the {\it generating distribution\/} of the convex
body~$K$. It is well known that positive distributions are in fact positive
measures. Therefore in the context of zonoids Weil's result is particularly
useful because it provides a priori a functional whose positiveness is to
be checked. For an illustration of this technique, see \cite{4}, Th.5.1.\par
\def\spc{C_e^{\8}(S^{n-1})}
The problem is now to prove that for the norm of the three dimensional
barrel zonoid, the generating distribution $T^{-1}(\norma{\cdot})$ is positive.
To this end an inversion formula for the cosine
transform shall be used, which involves
the {\it spherical Radon transform\/} $R:\spc\to\spc$, defined by
$$(Rf)(u)={1\over\omega_{n-1}}\int_{S^{n-2}\cap 
u^{\perp}}f(v)\,d\lambda_{n-2}(v),\quad u\in S^{n-1},$$
where $\omega_{n-1}$ is the total spherical Lebesgue measure of the unit
sphere in $\R^{n-1}$. Let $\Delta_n$ denote the spherical Laplace
operator on $S^{n-1}$. In \cite{4}, Goodey and Weil prove the following 
inversion
formula:
$$T^{-1}={1\over 2\omega_{n-1}}(\Delta_n+n-1)R^{-1}.\equano$$
It is well known that the Radon transform is a self--adjoint
continuous bijection of $\spc$ to itself. (see \cite{6}). Since~$T$
and $\Delta_n$ also have this
property, the inversion formula can be applied to the dual space
of its natural domain, that is, to the space of even distributions. In 
particular it can be applied to any given norm, restricted to the sphere.\par
As for the inversion of the Radon transform, it is explained by Gardner in 
\cite{2},
that if $f$ is a rotationally symmetric function
on $S^{n-1}$ and $f=Rg$ then~$g$
is also rotationally symmetric and the equation $f=Rg$ becomes
$$f(\sin^{-1}x)={2\omega_{n-2}\omega_{n-1}^{-1}\over 
x^{n-3}}\int_0^xg(\cos^{-1}t)(x^2-t^2)^{(n-4)/2}\,dt,\equano$$
for $0<x\leq 1$ and $g(\pi/2)=f(0)$. There is an inversion formula
for this equation (see \cite{2}). However, for $n=4$ it is trivial to invert the
equation $\numnow$ because if $xf(\sin^{-1}x)$ is differentiable, 
then $\numnow$ immediately
implies:
$$g(\cos^{-1}x)={d\over dx}(xf(\sin^{-1}x)),\equano$$
for $0\leq x\leq 1$. In proving the next claim, this formula
will be used.

\proclaim{Claim} The barrel zonoid $\Cal B_{4,r}$ is a polar of a zonoid if
and only if $r\leq 1$.\endproclaim

\demo{Proof}
$\numnow$ shall be now applied for
$n=4$ and for $f=f_r$, the norm of $\Cal B_{4,r}$, given by $\numback 5$.
Let $x_r=\pr$. Then
$$f_r(\sin^{-1}x)=\cases \sqrt{1-x^2},& {\text {if }} 0\leq x\leq x_r\cr
                       {1\over rx+\sqrt{1-r^2+r^2x^2}}, & {\text {if }} x_r\leq 
x\leq 1\endcases$$
The function $f_r(\sin^{-1}x)$ has a continuous derivative in $[0,1]$.
Hence, by $\numnow$,
$$g(\cos^{-1}x)=
\cases{1-2x^2\over\sqrt{1-x^2}},& {\text {if }} 0\leq x\leq x_r\cr
        {1-r^2\over A(x,r)(rx+A(x,r))^2},& {\text {if }} x_r\leq x\leq 
1\endcases$$
where $A(x,r)=\sqrt{1-r^2+r^2x^2}$. \par 
In cylindrical coordinates, $u=(\sqrt{1-x^2}\xi,x)$,
where $\xi\in S^2$, the $4$--dimensional spherical Laplacian
is given by 
$$\Delta_4=(1-x^2){\partial^2\over\partial x^2}-3x{\partial\over\partial 
x}+{1\over 1-x^2}\Delta_3.$$
In this formula $\Delta_3$ is applied to coordinates of~$\xi$; these are
independent of~$x$. Hence, when applying the Laplacian to a rotationally
symmetric function, the term containing $\Delta_3$ disappears. Therefore
the differential operator which
is to be applied to $g(\cos^{-1}x)$ is given by:
$$D={1\over 8\pi}\((1-x^2){d^2\over dx^2}-3x{d\over dx}+3\).$$
The calculation of the corresponding derivatives of $G(x)=g(\cos^{-1}x)$ has to 
be done in the distribution sense. The first derivative is an absolutely
continuous measure whose density is given by:
$${dG\over dx}=\cases{2x^3-3x\over(1-x^2)^{3/2}},& {\text {if }} 0\leq x<x_r\cr
                       {r(2A(x,r)+rx)(rx-A(x,r))\over A^3(x,r)(A(x,r)+rx)},& 
{\text {if }}
x_r<x\leq 1\endcases$$
Due to the jump at the point $x=x_r$, the second derivative is
a sum of a continuous measure and a measure concentrated at the point $x=x_r$.
(see \cite{3},\S2). Therefore,
$${d^2 G\over dx^2}=\cases{3\over(1-x^2)^{5/2}},& {\text {if }} 0\leq x<x_r\cr
                            {3r^2(1-r^2)\over A^5(x,r)},& {\text {if }} 
x_r<x\leq 1.\endcases\quad
+\quad c(r)\,\delta(x-x_r)$$
 The constant $c(r)$ is given by
$$c(r)=\lim_{x\to x_r^+}{dG\over dx}-\lim_{x\to x_r^{-}}{dG\over 
dx}=r(r^2+1)^2.$$
Having all the derivatives the differential
operator~$D$ can be applied directly. The result is
$$T^{-1}(f_r)=D(G)=\cases0,& {\text {if }} 0\leq x<x_r\cr
                         {3(1-r^2)\over 8\pi A^5(x,r)},& {\text {if }}
x_r<x\leq 1\endcases\quad+\quad{r(r^2+1)\over 8\pi}\delta(x-x_r).\equano$$

Evidently, the generating distribution is a positive measure if and only if 
$r\leq 1$.
The proof of the claim is complete.\qed\enddemo
\vfill\eject
{\bf Remarks}\par

{\bf 1.} In the special case where $n=4,r=1$, the measure which is obtained
in $\numnow$ is particularly simple. It is therefore easy to check directly that
it represents the polar of $\Cal B_4$. Here is the
calculation.\par
The formula $\numnow$ shows that the measure in question is concentrated on
the set $\{v\in S^3:|v_4|=\sqrt{1/2}\}$. By rotational symmetry, its restriction
to each one of the $3$--dimensional spheres comprising its support is
proportional to the corresponding Lebesgue spherical measure. In order to
check this, let $e_4=(0,0,0,1)$ and consider the integral
$$\int\limits_{S^3\cap e_4^{\perp}}|\innerprod 
{y+e_4}{u}|\,d\lambda_2(y).\equano$$
Using invariance the point $u$ can be replaced by the point
$(0,0,\sqrt{1-t^2},t)$, where ${t=u_4}$. Applying spherical coordinates to 
$\numnow$
yields the following integral,
$$\aligned &2\pi\int_0^{\pi}|\sqrt{1-t^2}\cos\psi+t|\sin\psi
\,d\psi\\
&=2\pi{\cases
2|t|,            &{\text {if }} |t|\geq {1/\sqrt{2}}\\
{1\over\sqrt{1-t^2}},&{\text {if }} |t|\leq 1/\sqrt{2}.\endcases}\cr
&=4\pi f_1(\cos^{-1}|t|),\endaligned$$
where $f_1$ is the norm of the barrel given by (2.1). 
Hence a (positive) multiple of the spherical Lebesgue
measure on each one of the spheres $\{u\in S^3:u_4=\pm\,\sqrt{1/2}\}$, 
represents the polar
$\Cal B_4^{\circ}$, as was required to check.\medskip
{\bf 2.} Since $\Cal B_{n,r}$ is a central section of $\Cal B_{n+1,r}$ on which
there exists an orthogonal projection, the fact that $\Cal B_{n+1,r}^{\circ}$
is a zonoid implies the same for $\Cal B_{n,r}^{\circ}$. 
Consequently, $\Cal B_{3,r}^{\circ}$
is a zonoid for $r\leq 1$. For the special case $r=1$ it is possible to apply
the same reasoning as above and obtain an explicit formula for the 
generating distribution of $\Cal B_3^{\circ}$. Its density is
given by
$$T^{-1}f(\varphi)=\cases c{\cos^2\varphi+\sqrt{\cos 2\varphi}
\over (1+\sqrt{\cos 2\varphi})^2\cos^3\varphi\sqrt{\cos 2\varphi}},
& {\text {if }} 0\leq \varphi<{\pi\over 4}\cr
0,& {\text {if }} {\pi\over 4}\leq \varphi\leq {\pi\over 2}\endcases$$
Here $c>0$ is a positive constant. Hence the generating distribution of $\Cal 
B_3^{\circ}$
has an $L_1$ (but not $L_2$) density.\medskip
{\bf 3.} For $n\geq 6$ and $r>0$, the polar of $\Cal B_{n,r}$ is not a
zonoid. Indeed, the calculation of this section for
the case $n=6$, results in a generating distribution that
involves a derivative of a measure concentrated at a point. Consequently,
the generating distribution of $\Cal B_{6,r}^{\circ}$ is not a measure.
A generalization of this calculation to higher dimensions can be
used to answer a question raised by Goodey and Weil. For more details see 
\cite{9}.\par
It is plausible that non-smooth zonoids whose polars are zonoids exist
in every dimension. However, the author is not aware of any examples other
than the barrel zonoids.\bigskip
I thank Professor J. Lindenstrauss for useful discussions.
\vfill\eject

\enddocument